\def\A{\leavevmode\setbox0\hbox{A}\lower1.4ex\hbox to\wd0{\hss`}\kern-.9\wd0A}
\def\E{\leavevmode\setbox0\hbox{E}\lower1.4ex\hbox to\wd0{\hss`\/}\kern-.9\wd0E}
\def\a{\leavevmode\setbox0\hbox{a}\lower1.4ex\hbox to\wd0{\hss`\/}\kern-\wd0a}
\def\e{\leavevmode\setbox0\hbox{e}\lower1.4ex\hbox to\wd0{\hss`\/}\kern-\wd0e}
\newcommand{\be}{\begin{equation}}
\newcommand{\ee}{\end{equation}}
\newcommand{\ba}{\begin{array}}
\newcommand{\ea}{\end{array}}
\newcommand{\beqn}{\begin{eqnarray}}
\newcommand{\eeqn}{\end{eqnarray}}

\newcommand{\zero}{\setcounter{equation}{0} \par}

\def\znakr{\raise1.5pt\hbox{\symb\char66\kern-2pt\char74}}
\def\znakl{\raise1.5pt\hbox{\symb\char73\kern-2pt\char67}}

\def\normalsize{
\setlength{\textheight}{23cm}
\setlength{\textwidth}{15cm}
\setlength{\topmargin}{-2.0cm}
\setlength{\hoffset}{-0.5cm}
\setlength{\leftmargin}{-1cm}
\setlength{\rightmargin}{2.0cm}}

\documentstyle{article}
\normalsize
\begin{document}
%%%%%%%%%%%%%%%%%%%%%%%%%%%%%%%%%%%%%%%%%%%%%%%%%%%%%%%%%%%%
%\baselineskip=24pt
%\baselineskip=20pt
%%%%%%%%%%%%%%%%%%%%%%%%%%%%%%%%%%%%%%%%%%%%%%%%%%%%%%%%%%%%
%%%%%%%%%%%%%%%%%%%%%%%%%%%%%%%%%%%%%%%%%%%%%%%%%%%%%%%%%%%%
%%%%%%%%%%%%%%%%%%%%%%%%%%%%%%%%%%%%%%%%%%%%%%%%%%%%%%%%%%%%
%%%%%%%%%% BEGIN TITLE PAGE %%%%%%%%%%%%%%%%%%%%%%%%%%%%%%%%
%%%%%%%%%%%%%%%%%%%%%%%%%%%%%%%%%%%%%%%%%%%%%%%%%%%%%%%%%%%%
\title{The infinitesimal form of induced representation of
the $\kappa -$Poincar\'e group}
\author{ Karol Przanowski\thanks{Supported by
\L{}\'od\'z University grant no 487} \\
Department of Field Theory \\
University of \L \'od\'z \\
ul.Pomorska 149/153, 90-236 \L \'od\'z , Poland}
%%%%%%%%%%%%%%%%%%%%%%%%%%%%%%%%%%%%%%%%%%%%%%%%%%%%%%%%%%%%
\date{}
\maketitle
\setcounter{section}{0}
\setcounter{page}{1}
%%%%%%%%%%%%%%%%%%%%%%%%%%%%%%%%%%%%%%%%%%%%%%%%%%%%%%%%%%%%
\begin{abstract}
The infinitesimal form of the induced representation of
the $\kappa -$Poincar\'e group is constructed. The infinitesimal
action of the $\kappa -$Poincar\'e group on the $\kappa -$Minkowski space
is described. The actions of these two infinitesimal forms on
the solution of Klein-Gordon equation are compared.
\end{abstract}
%%%%%%%%%%%%%%%%%%%%%%%%%%%%%%%%%%%%%%%%%%%%%%%%%%%%%%%%%%%%
%%%%%%%%%%% END TITLE PAGE %%%%%%%%%%%%%%%%%%%%%%%%%%%%%%%%%
%%%%%%%%%%%%%%%%%%%%%%%%%%%%%%%%%%%%%%%%%%%%%%%%%%%%%%%%%%%%
%%%%%%%%%%%%%%%%%%%%%%%%%%%%%%%%%%%%%%%%%%%%%%%%%%%%%%%%%%%%
%%%%%%%%%%%%%%%%%%%%%%%%%%%%%%%%%%%%%%%%%%%%%%%%%%%%%%%%%%%%

\section{Introduction}
\zero

Recently, considerable interest has been paid to the deformations of group
and algebras of space-time symmetries~\cite{w1}. An interesting deformation of
the Poincar\`e algebra~\cite{w2},~\cite{defwayla} as well
as group~\cite{zak} has been introduced which
depend on the dimensional deformation parameter $\kappa $; the relevant
objects are called  $\kappa -$Poincar\`e algebra and $\kappa -$Poincar\`e
group, respectively. Their structure was studied in some detail and many
of their properties are now well understood. In particual, the induced
representations of the deformed group were found~\cite{repre} and the duality
between $\kappa -$Poincar\`e group and $\kappa -$Poincar\`e algebra
was also given~\cite{w11}. Having the representations of the $\kappa -$Poincar\`e
group and duality relations one can consider the infinitesimal form of
the representation; in order to check whether we obtain the representation
of the $\kappa -$Poincar\`e algebra. This is nontrivial as it is known
from the construction of the induced representations the action of the
Lorentz group on $q-$space is standard. Therefore the support spaces
are described by the classical equation $q^2 = const$. On the other hand
$q^2$ is not the Casimir operator of $\kappa -$Poincar\`e algebra.

In section~\ref{r1} we describe the general definition of the induced representation
of quantum group,
than we find the infinitesimal form of the induced representation
 of $\kappa -$Poincar\'e group constructed by~\cite{repre}.
We show that in the massive case the infinitesimal
form of the induced representation is the representation
of the $\kappa -$Poincar\`e algebra.
Next in section~\ref{r2} we consider the infinitesimal action
of ${\cal P}_\kappa $ on
$\kappa -$Minkowski space ${\cal M}_\kappa $,
and describe the Klein-Gordon equation.
And on the end we compare actions of our two infinitesimal forms
over the solution K-G equation.

In that paper we asume that $g_{\mu \nu }$ is diagonal $(+,-,-,-)$
tensor matric.

Let us remind the definition of $\kappa -$Poincar\`e group.
The $\kappa -$Poincar\`e group ${\cal P}_\kappa $ is the Hopf *-algebra
generated by selfadjoint elements $\Lambda ^\mu _{\ \nu },\ v^\mu $
subject to the following relations:
\beqn
\lbrack  \Lambda ^\alpha _{\ \beta },v^\varrho  \rbrack &=& - {i\over \kappa }
((\Lambda ^\alpha _{\ 0}-\delta ^\alpha _{\ 0})
\Lambda ^\varrho _{\ \beta }+(\Lambda _{0\beta }-
 g_{0\beta })g^{\alpha \varrho }), \nonumber \\
\lbrack  v^\varrho ,v^\sigma  \rbrack  &=& {i\over \kappa }
(\delta ^\varrho _{\ 0}v^\sigma  -
\delta ^\sigma _{\ 0}v^\varrho ), \nonumber \\
\lbrack  \Lambda ^\alpha _{\ \beta }, \Lambda ^\mu _{\ \nu } \rbrack
&=& 0.  \nonumber
\eeqn

The comultiplication, antipode and counit are defined as follows:
\beqn
\Delta \Lambda ^\mu _{\ \nu } &=& \Lambda ^\mu _{\ \alpha } \otimes
\Lambda ^\alpha _{\ \nu },\nonumber \\
\Delta v^\mu  &=& \Lambda ^\mu _{\ \nu }\otimes v^\nu  +
v^\mu \otimes I, \nonumber \\
S(\Lambda ^\mu _{\ \nu }) &=& \Lambda ^{\ \mu }_\nu ,\nonumber \\
S(v^\mu )&=& -  \Lambda ^{\ \mu }_\nu  v^\nu ,\nonumber \\
\varepsilon (\Lambda ^\mu _{\ \nu })&=& \delta ^\mu _{\ \nu },\nonumber \\
\varepsilon (v^\mu )&=& 0  . \nonumber
\eeqn

Its dual structure, the $\kappa -$Poincar\`e algebra $\tilde {\cal P}_\kappa $
(in the Majid and Ruegg basis~\cite{w9}) is a quantized universal envoloping
algebra  in the sense of Drinfeld~\cite{drin} described by the
following relations:
\beqn
\lbrack  M_{ij},P_0 \rbrack  &=& 0 ,\nonumber \\
\lbrack M_{ij},P_k  \rbrack  &=& i(g_{jk}P_i - g_{ik}P_j), \nonumber \\
\lbrack M_{i0},P_0  \rbrack  &=& iP_i ,\nonumber \\
\lbrack M_{i0},P_k  \rbrack  &=& -i{\kappa \over 2}
g_{ik}(1-e^{{-2\over \kappa }P_0})+
{i\over 2\kappa }g_{ik}P^r P_r -
{i\over \kappa } P_i P_k ,\nonumber \\
\lbrack P_\mu ,P_\nu   \rbrack  &=& 0 ,\nonumber \\
\lbrack M_{i j },M_{r s }  \rbrack
&=& i(g_{i s }M_{j r }-g_{j s }M_{i r }
+g_{j r }M_{i s }-g_{i r }M_{j s }) ,\nonumber \\
\lbrack M_{i0}, M_{rs} \rbrack &=& -i(g_{is}M_{r0}
- g_{ir}M_{s0}) ,\nonumber  \\
\lbrack M_{i0}, M_{j0} \rbrack &=& -iM_{ij} .\label{komut}
\eeqn
The coproducts, counit and antipode:
\beqn
\Delta P_0 &=& I\otimes P_0 + P_0 \otimes I ,\nonumber \\
\Delta P_k &=& P_k \otimes e^{-{P_0\over \kappa }}+
I\otimes P_k ,\nonumber \\
\Delta M_{ij} &=& M_{ij}\otimes I + I\otimes M_{ij} ,\nonumber \\
\Delta M_{i0} &=& I\otimes M_{i0} + M_{i0}\otimes
e^{-{P_0\over \kappa }}+ {1\over \kappa }M_{ij}\otimes P_j ,\nonumber \\
\varepsilon (M_{\mu \nu }) &=& 0; \ \ \
\varepsilon (P_\nu ) = 0, \nonumber \\
S(P_0) &=& -P_0 ,\nonumber \\
S(P_i) &=& -e^{P_0\over \kappa }P_i ,\nonumber \\
S(M_{ij}) &=& -M_{ij} ,\nonumber \\
S(M_{i0}) &=& -e^{P_0\over \kappa }(M_{i0}-
{1\over \kappa }M_{ij}P_j) ,\nonumber
\eeqn
where $i,j,k,r,s = 1,2,3$.

Fact that the $\kappa -$Poincar\`e algebra is dual to the 
$\kappa -$Poincar\`e group was proved by~\cite{w11}
The fundamental duality relations read:
\beqn
< P_\mu , f(v) > &=&  i \left. {\partial \over \partial v^\mu } f(v) \right|_{v=0}
\nonumber  \\
< M_{\mu \nu } , f(\Lambda ) > &=& i \left. ({\partial \over \partial
\Lambda ^{\mu \nu }}
-{\partial \over \partial \Lambda ^{\nu \mu }}) f(\Lambda ) \right|_{\Lambda =I}
%\nonumber
\label{dual}
\eeqn

\section{The infinitesimal form of the induced representation \label{r1}}
\zero

Let us recall the definition of the representation of a quantum group
${\cal A}(G)$, acting in the linear space ${\cal V}$. It is simply a map
\beqn
\varrho : {\cal V} \to {\cal V} \otimes {\cal A}(G)
\nonumber
%\label{defro}
\eeqn
satisfying
\beqn
(I\otimes \Delta )\otimes \varrho = (\varrho   \otimes I)\otimes \varrho
\nonumber
%\label{lacz}
\eeqn
The induced representations are defined as follow~\cite{repre},~\cite{ruiz}:
given any quantum group
${\cal A}(G)$, its quantum subgroup ${\cal A}(H)$ and the representation
$\varrho $ of the latter acting in the linear space ${\cal V}$, we consider
the subspace of the space ${\cal V}\otimes {\cal A}(G)$ defined by the
coequivariance condition:
$$
\tilde {\cal V} = \left\{ F \in {\cal V}\otimes {\cal A}(G):
id\otimes (\Pi \otimes id)\circ \Delta _G F =
(\varrho \circ id)F  \right\}
$$
where $\Pi : {\cal A}(G) \to {\cal A}(H)$ is epimorfizm defining the subgroup
${\cal A}(H)$.

The induced representation is defined as a right action:
\beqn
\tilde \varrho &:& \tilde {\cal V} \to \tilde {\cal V} \otimes {\cal A}(G) \nonumber \\
\tilde \varrho  &=& id \otimes \Delta
\nonumber
%\label{induk}
\eeqn
In the paper~\cite{repre} Ma\'slanka obtained the following form
of the induced representation in the massive case: \par
{\it (elements of $\tilde {\cal V} $ we
write $f(q_\mu ) = e_i \otimes f_i (q_\mu ) $,
where $\{ e_i \}$ are a basis of ${\cal V} $)}
\newcommand{\sh}{\sinh({m\over \kappa })}
\newcommand{\ch}{\cosh({m\over \kappa })}
\beqn
\varrho _{\cal R} : f_i(q_\nu ) &\longmapsto & {\cal D}_{ij}({\cal R}
(\widetilde{q},\Lambda ))\cdot \label{ror} \\
&\ & \exp (-i\kappa \ln (\ch - {q_0\over m}\sh )\otimes v^0)\cdot  \nonumber \\
&\ & \exp ({i\kappa \sh q_k\over m\ch -q_0 \sh }\otimes v^k)
f_j (q_\mu \otimes \Lambda ^\mu _{\ \nu }) \nonumber
\eeqn
where:
\beqn
q_\mu &=& m\Lambda ^0_{\ \mu } \nonumber  \\
\widetilde{q_0} &=& {m q_0 \ch -m^2 \sh \over m \ch - q_0 \sh }
\nonumber  \\
%\label{pzero} \\
\widetilde{q_k} &=& {m q_k \over m \ch - q_0 \sh }
\nonumber
%\label{pk}
\eeqn
and $f(q)$ are the square integrable functions defined on the hyperboloid
$q^2=m^2$, $( q_0 = \sqrt{q_i q_i + m^2} )$ and taking values in the vector space
carrying the unitary
representation of the rotation group, the matrices ${\cal D}_{ij}$ are
constructed in the same way as the matrices of the representation of
classical orthogonal group and ${\cal R}(\widetilde{q},\Lambda )$ is
a classical
Wigner rotation corresponding to the momentum $\widetilde{q}$ and transformation
$\Lambda $,~\cite{repre}.
Of course, the right-hand side of eq.(\ref{ror}) is to be understood here
as an element of the tensor product of the algebra of function on the
hyperboloid $q^2=m^2$ and the group ${\cal P}_\kappa $.
Following Woronowicz~\cite{sudwa} we define
the infinitesimal form of the induced representation.
For any element $X$ of enveloping algebra
$\tilde {\cal P}_\kappa $ and for any $f \in {\cal P}_\kappa $ if
$$
X(f)=<X,f>
$$
and
$$
\varrho_{\cal R} : {\cal V} \to {\cal V} \otimes {\cal P}_\kappa
$$
we define
$$
\widetilde{X} : {\cal V} \to {\cal V}
$$
by
$$
\widetilde{X} = (I\otimes X)\circ \varrho_{\cal R}
$$

Using the duality relations (\ref{dual}) after some calculi we arrive
at the following formulas describing the infinitesimal form of
our representation:
\beqn
\widetilde{M}_{ij} &=& i(q_i {\partial \over \partial q^j}-
q_j {\partial \over \partial q^i})+\varepsilon _{ijk}s_k \nonumber \\
\widetilde{M}_{i0} &=& -i q_0 {\partial \over \partial q^i}
+\varepsilon _{ijk}{q_j s_k \over q_0 + m} \nonumber \\
\widetilde{P}_0 &=& p_0 = \kappa \ln (\ch - {q_0\over m}\sh ) \nonumber \\
\widetilde{P}_j &=& p_j= -\kappa {\sh q_j\over m \ch - q_0 \sh }
\label{infinit1}
\eeqn
where $s_k,\ (k=1,2,3) $ are the infinitesimal forms of
representation ${\cal D}_{ij}$ ( The representation of $SU(2)$ algebra ).

It is easy to check that our operators satysfying relations of
the $\kappa -$Poincar\'e algebra (\ref{komut}), so the infinitesimal
form is the representation of the $\kappa -$Poincar\'e algebra.

\section{The $\kappa -$Minkowski space, the infinitesimal action, K-G equation
\label{r2}}
\zero

The $\kappa -$Minkowski space~(\cite{zak},~\cite{note})  ${\cal M}_\kappa $
is a universal $*$-algebra
with unity generated by four selfadjoint elements $x^\mu$ subject to the
following conditions:
$$
\lbrack x^\mu, x^\nu \rbrack  =  {i\over \kappa} (\delta ^\mu_0 x^\nu
- \delta ^\nu_0 x^\mu ) .
$$
Equipped with the standard coproduct:
$$
\Delta  x^\mu =  x^\mu \otimes  I + I \otimes  x^\mu ,
$$
antipode $S(x^\mu) = - x^\mu$ and counit $\varepsilon (x^\mu) = 0$
it becomes a quantum group.

The product of generators $x^\mu$ will be called normally ordered if all
$x^0$ factors stand leftmost. This definition can be used to ascribe a
unique element $:f(x) :$ of ${\cal M}_\kappa $ to any polynomial
function of four
variables $f$. Formally, it can be extended to any analytic function $f$.

Let us now define the infinitesimal action of ${\cal P}_\kappa $
on ${\cal M}_\kappa $.
The $\kappa -$Minkowski space carries a left-covariant action of
$\kappa -$Poincar\'e group ${\cal P}_\kappa $,\
$\varrho _L : {\cal M}_\kappa \to {\cal P}_\kappa \otimes {\cal M}_\kappa $,
given by
\be
\varrho _L (x^\mu ) = \Lambda ^\mu _{\ \nu } \otimes x^\nu +v^\mu \otimes I.
\label{rorl}
\ee
Let $X$ be
any element of the Hopf algebra dual to ${\cal P}_\kappa $ -- the
$\kappa -$Poincar\'e algebra $\tilde {\cal P}_\kappa $.
The corresponding infinitesimal action:
$$
\widehat{X} : {\cal M}_\kappa  \to {\cal M}_\kappa
$$
is defined as follows: for any $f \in {\cal M}_\kappa $,
$$
\widehat{X} f = (X \otimes  I) \circ \varrho _L (f).
$$
The following forms of generators were obtained in~\cite{note}:
\beqn
\widehat{P}_\mu  :f:\ &=& \
: i {\partial f\over \partial x^\mu } :  \nonumber  \\
\widehat{M}_{ij} :f: \ &=& \ : -i
(x_i {\partial \over \partial x^j}
-x_j {\partial \over \partial x^i})f :   \nonumber  \\
\widehat{M}_{i0} :f: \  &=&
\ :\left[ ix^0 {\partial \over \partial x^i} - x_i \left(
 {\kappa \over 2}(1-e^{-{2i\over \kappa }{\partial \over \partial x^0}})
- {1\over 2\kappa }\Delta   \right) \right. \nonumber \\
&+& \left. {1\over \kappa }x^k {\partial ^2\over \partial x^k \partial x^i}
 \right] f: \nonumber
\eeqn
In that case these operators are not satysfying relation of
$\kappa -$Poincar\'e algebra (\ref{komut}),
becouse the action of $\kappa -$Poincar\`e algebra (group) on the
$\kappa -$Minkowski space is antirepresentation (not representation).
It is clear from the following equation:
$$
< AB f(P), \varphi (x) > = <f(P), \widehat{B}\widehat{A} \varphi (x)>
$$
for $A,B,P \in \tilde {\cal P}_\kappa $.

The deformed Klein-Gordon equation we write in the two equivalent forms:
$$
\left( \partial +{m^2\over 8} \right) f = 0,
$$
or
$$
\left[ \partial ^2_0 - \partial ^2_i + m^2 (1+{m^2\over 4\kappa ^2})
\right] f = 0,
$$
where operators $\partial _0 ,\ \partial _i,\ \partial $ are defined:
\beqn
\partial _0 :f:\ &=& \ :\left(  \kappa \sin ({1\over \kappa }
{\partial \over \partial x^0})+{i\over 2\kappa }
e^{{i\over \kappa }{\partial \over \partial x^0}}\Delta
\right) f: \nonumber  \\
\partial _i :f: \ &=& \ :(e^{{i\over \kappa }{\partial \over \partial x^0}}
{\partial \over \partial x^i}) f: \nonumber \\
\partial  :f:\ &=& \ : \left( {\kappa ^2\over 4}
(1-\cos ({1\over \kappa }\partial _0)) - {1\over 8}e^{{i\over \kappa }
\partial _0}\Delta  \right) f:\ . \nonumber
\eeqn

We can write the solution of K-G equation
in the following wave function~\cite{luk}:
\beqn
\Phi (x^\mu ) &=& \ : \int {d^3 \vec q\over q^0} a(\vec q)
e^{-ip_\mu (q)x^\mu } \ : \nonumber
\eeqn
where $p_\mu $ is deformed of $q_\mu $ defined in eq.(\ref{infinit1}).

For $\widehat{X}=
\widehat{M}_{i0},\ \widehat{M}_{ij},\ \widehat{P}_\mu $, 
following the paper~\cite{luk},
let us define the operators ${\cal X}(q)$ by the following relation:
\beqn
\widehat{X}(q_\nu ) \Phi (x^\mu ) &=& \ : \int {d^3 \vec q\over q^0}
\left\{ {\cal X}(q_\nu ) a(\vec q) \right\}
e^{-ip_\mu (q)x^\mu } \ : \nonumber
%\label{cal}
\eeqn

It is easy to see that:
\beqn
{\cal M}_{ij}(q) &=& -i(q_i {\partial \over \partial q^j}-
q_j {\partial \over \partial q^i}) = - \widetilde{M}_{ij}(q) \nonumber \\
{\cal M}_{i0}(q) &=& i q_0 {\partial \over \partial q^i} =
- \widetilde{M}_{i0}(q) \nonumber \\
{\cal P}_\mu (q) &=& p_\mu (q)= \widetilde{P}_\mu (q) \nonumber
\eeqn
The last relations compare generators defined in~\cite{luk} from
with our found from induced representation.
\section{Acknowledgment}
\zero

The numerous discussions with prof.J.Lukierski and 
dr Pawe\l{}  Ma\`slanka are kindly acknowledged.

%%%%%%%%%%%%%%%%%%%%%%%%%%%%%%%%%%%%%%%%%%%%%%%%%%%%%%%%


\begin{thebibliography}{99}
%%%%%%%%%%%%%%%%%%%%%%%%%%%%%%%%%%%%%%%%%%%%%%%%%%%%%%%%
\newcommand{\byauthors}[1]{#1 }
\newcommand{\journal}[1]{ #1 }
\newcommand{\reftitle}[1]{{\it #1} }
\newcommand{\volumin}[1]{{\bf #1} }
\newcommand{\eref}{.}
\newcommand{\refyear}[1]{(#1), }
%%%%%%%%%%%%%%%%%%%%%%%%%%%%%%%%%%%%%%%%%%%%%%%%%%%%%%%%
%\bibitem{}\byauthors{}\reftitle{}\journal{}\refyear{}\volumin{}\eref
%%%%%%%%%%%%%%%%%%%%%%%%%%%%%%%%%%%%%%%%%%%%%%%%%%%%%%%%

\bibitem{zak}\byauthors{S.Zakrzewski}\reftitle{}\journal{J.Phys.}\refyear{1994}
\volumin{A}27\eref

\bibitem{note}\byauthors{S.Giller, C.Gonera, P.Kosi\'nski, P.Ma\'slanka}
\reftitle{A note on geometry of $\kappa -$Minkowski space}
- in print \journal{Acta. Phys. Pol.}\volumin{B}\eref

\bibitem{repre}\byauthors{P.Ma\'slanka}
\reftitle{The induced representations of the $\kappa -$Poincar\'e group.}
\journal{J. Math. Phys.}35(4)\refyear{1994}5047-5056\eref

\bibitem{w11}\byauthors{P.Kosi\'nski,  P.Ma\'slanka
}\reftitle{preprint IMU\L \/
3\slash 94}or \volumin{Q-ALG}9411033\eref

\bibitem{luk}\byauthors{P.Kosi\'nski, J.Lukierski, P.Ma\'slanka, A.Sitarz}
- paper in preparation\eref

\bibitem{defwayla}\byauthors{P.Kosi\'nski, P.Ma\'slanka}
\reftitle{The $\kappa -$Weyl group and its algebra}
\journal{in "From Field Theory to Quantum Groups"}\volumin. on 60
anniversary of J.Lukierski, World Scientific Singapur, 1996
or \volumin{Q-ALG}9512018\eref


\bibitem{w1}\byauthors{W.B.Schmidke, J.Weiss, B.Zumino}
\journal{Zeitschr. f. Physik}
\volumin{52}\refyear{1991}472\eref

\byauthors{U.Carow-Watamura, M.Schliecker, M.Scholl, S.Watamura}
\journal{Int. J. Mod. Phys.}\volumin{A\,6}\refyear{1991}3081\eref

\byauthors{S.L. Woronowicz}\journal{Comm. Math. Phys.}
\volumin{136}\refyear{1991}399\eref

\byauthors{O. Ogievetsky, W.B. Schmidke, J. Weiss, B. Zumino}
\journal{Comm. Math.  Phys.}\volumin{150}\refyear{1992}495\eref

\byauthors{M. Chaichian, A.P. Demichev}\journal{Proceedings of the
Workshop: "Generalized symmetries in Physics", Clausthal}\refyear{1993}\eref

\byauthors{V. Dobrev}\journal{J.  Phys.}
\volumin{A\,26}\refyear{1993}1317\eref

\byauthors{L. Castellani}\journal{in "Quantum Groups" Proceedings of XXX
Karpacz Winter School of Theoretical Physics, Karpacz 1994, PWN 1995, p.
13}\eref


\bibitem{w2}\byauthors{J. Lukierski, A. Nowicki, H. Ruegg, V. Tolstoy}
\journal{Phys. Lett.}\volumin{B\,264}\refyear{1991}331\eref

\byauthors{J. Lukierski, A. Nowicki, H. Ruegg}\journal{Phys.
Lett.}\volumin{B\,293}\refyear{1993}344\eref

\byauthors{S. Giller,  P. Kosi\'nski, J. Kunz,  M. Majewski, P.
Ma\'slanka}\journal{Phys. Lett.}\volumin{B\,286}\refyear{1992}57\eref


\bibitem{w9}\byauthors{S. Majid, H. Ruegg}\journal{Phys. Lett.}
\volumin{B\,334}
\refyear{1994}348\eref

\byauthors{Ph. Zaugg}\journal{preprint MIT--CTP--2353}\refyear{1994}\eref


\bibitem{drin}\byauthors{V.G.Drinfeld}
\journal{Proc. Int. Congr. of Math.}Berkeley\refyear{1986}p.798\eref


\bibitem{ruiz}\byauthors{A.Gonzales-Ruiz, L.A.Ibort}
\journal{Phys. Lett.}\refyear{1992}\volumin{B 296}104\eref


\bibitem{sudwa}\byauthors{S.L.Woronowicz}
\journal{Publ. RIMS, Kyoto Univ.}\refyear{1987}\volumin{23}117-181\eref


\end{thebibliography}
\end{document}